\documentclass[runningheads]{llncs}
\usepackage[T1]{fontenc}
\usepackage{graphicx}
\usepackage{url}
\usepackage{amsmath}    % basic ams math environments and symbols
\usepackage{amssymb}    % ams symbols and theorems
\usepackage{physics}
\usepackage{xcolor}
\usepackage{enumitem}

\renewcommand{\L}{\mathcal{L}}
\renewcommand{\b}{\boldsymbol}

\newcommand{\x}{\b{x}}
\newcommand{\w}{\b{w}}

\begin{document}
\title{Spatially-varying meshless approximation method for enhanced computational efficiency}
\titlerunning{Spatially-varying meshless approximation method}
\author{Mitja Jan\v{c}i\v{c}\inst{1,2}\orcidID{0000-0003-1850-412X} \and
	Miha Rot\inst{1,2}\orcidID{0000-0002-8869-4507} \and
	Gregor Kosec\inst{1}\orcidID{0000-0002-6381-9078}}
\authorrunning{M.~Jan\v{c}i\v{c} et al.}
\institute{
	Jo\v{z}ef Stefan Institute, Parallel and Distributed Systems Laboratory, Ljubljana, Slovenia
	\email{\{mitja.jancic,miha.rot,gregor.kosec\}@ijs.si}
	\and
	Jo\v{z}ef Stefan International Postgraduate School, Ljubljana, Slovenia\\
}
\maketitle              % typeset the header of the contribution

\begin{abstract}
	In this paper, we address a way to reduce the total computational cost of meshless approximation by reducing the required stencil size through spatially varying computational node regularity. Rather than covering the entire domain with scattered nodes, only regions with geometric details are covered with scattered nodes, while the rest of the domain is discretized with regular nodes. A simpler approximation using solely monomial basis can be used in regions covered by regular nodes, effectively reducing the required stencil size and computational cost compared to the approximation on scattered nodes where a set of polyharmonic splines is added to ensure convergent behaviour.
	%Consequently, in regions covered with regular nodes the approximation using solely monomial basis can be performed, effectively reducing the required stencil size compared to the approximation on scattered nodes where a set of polyharmonic splines is added to ensure convergent behaviour.
	%In this paper, we propose to improve computational efficiency by spatially varying the meshless approximation methods. In particular, we propose to employ the computationally efficient collocation method on uniformly positioned discretization nodes and the more expensive but less sensitive to non-optimal node placement, radial basis function-generated finite differences (RBF-FD) on scattered nodes. 

	The performance of the proposed hybrid scattered-regular approximation approach, in terms of computational efficiency and accuracy of the numerical solution, is studied on natural convection driven fluid flow problems. We start with the solution of the de~Vahl~Davis benchmark case, defined on a square domain, and continue with two- and three-dimensional irregularly shaped domains. We show that the spatial variation of the two approximation methods can significantly reduce the computational demands, with only a minor impact on the accuracy.
\end{abstract}

\section{Introduction}
Although the meshless methods are formulated without any restrictions regarding the node layouts, it is generally accepted that quasi-uniformly-spaced node sets improve the stability of meshless methods~\cite{wendland2004scattered,liu2009meshfree}. Nevertheless, even with quasi-uniform nodes generated with recently proposed node positioning algorithms~\cite{slak2019generation,shankar2018robust,van2021fast}, a sufficiently large stencil size is required for stable approximation. A stencil with $n = 2\binom{m+d}{m}$ nodes is recommended~\cite{bayona2017role} for the local Radial Basis Function-generated Finite differences (RBF-FD)~\cite{tolstykh2003using} method in a $d$-dimensional domain for approximation order $m$. The performance of RBF-FD method --- with approximation basis consisting of Polyharmonic splines (PHS) and monomial augmentation with up to and including monomials of degree $m$ --- has been demonstrated with scattered nodes on several applications~\cite{slak2019adaptive,fornberg2015primer,zamolo2019solution}. On the other hand, approximation on regular nodes can be performed with considerably smaller stencil~\cite{kosec2018local} ($n=5$ in two-dimensional domain) using only monomial basis.

%In Radial Basis Function-generated Finite differences (RBF-FD)~\cite{tolstykh2003using}, a method that has been demonstrated on several applications using scattered nodes~\cite{slak2019adaptive,fornberg2015primer,toth2023h,zamolo2019solution}, a set of Polyharmonic splines (PHS) augmented with monomials up to certain order $m$ are used for approximation on $n = 2\binom{m+d}{m}$ stencil nodes in $d$-dimensional domain.

%Within the meshless community, a given domain is commonly discretized using only regular or only scattered nodes. Generally speaking, regular nodes are preferred because the computational efficiency of approximation methods compatible with regular nodes can be notably better than that of methods robust enough to support scattered nodes~\cite{patel2020meshless,liu2009meshfree}. However, discretizing the entire domain with only regular nodes can lead to significant discretization-related errors in presence of irregularly shaped domain boundaries. On the other hand, using only scattered nodes often turns out to be suboptimal from a computational efficiency point of view if the present irregularities are local. 

Therefore, a possible way to enhance the overall computational efficiency and consider the discretization-related error is to use regular nodes far away from any geometric irregularities in the domain and scattered nodes in their vicinity. A similar approach, where the approximation method is spatially varied, has already been introduced, e.g., a hybrid FEM-meshless method~\cite{elkadmiri2022hybrid} has been proposed to overcome the issues regarding the unstable Neumann boundary conditions in the context of meshless approximation. Moreover, the authors of~\cite{ding2004simulation,bourantas2018hybrid} proposed a hybrid of Finite Difference (FD) method employed on conventional cartesian grid combined with meshless approximation on scattered nodes. These hybrid approaches are computationally very efficient, however, additional implementation-related burden is required on the transition from cartesian to scattered nodes~\cite{javed2013hybrid}, contrary to the objective of this paper relying solely on the framework of meshless methods.

%could improve the overall computational efficiency while ensuring a satisfactory local field description. %Additionally, scattered nodes also allow to employ the $h$-refined solution procedures~\cite{slak2020adaptive} to further reduce the discretization-related errors near irregularly shaped boundaries.
%The generalized finite difference approach to meshless discretization offers a wide range of methods~\cite{slak2019refinedCauchy,tolstykh2003using,prax1996diffuse}, which differ mainly in the type and number of basis functions used. In this paper, we chose the two extremal second-order cases. The collocation method~\cite{kosec2018local,slak2019refinedCauchy}, which is computationally efficient but requires a uniform discretization, and the radial basis function-generated finite difference method (RBF-FD)~\cite{tolstykh2003using}, which is computationally expensive but stable on scattered nodes and resistant to non-ideal node configurations~\cite{bayona2017role,bayona2019insight,jancic2022stability}.

%We propose to spatially vary the discretization type and, accordingly, the approximation method. The collocation method is employed on regular nodes positioned sufficiently far away from irregularly shaped domain boundaries and RBF-FD on scattered nodes in their vicinity. 
In this paper we experiment with such \emph{hybrid scattered-regular} method with spatially variable stencil size on solution of natural convection driven fluid flow cases. The solution procedure is first verified on the reference de~Vahl~Davis case, followed by a demonstration on two- and three-dimensional irregular domains. We show that spatially varying the approximation method can have positive effects on the computational efficiency while maintaining the accuracy of the numerical solution.

\section{Numerical treatment of partial differential equations}
\label{sec:method}
%Generally speaking, obtaining a numerical solution to the governing system of partial differential equations (PDEs) consists of four steps. In the first step, the $d-$dimensional computational domain $\Omega\in \R^d$ is discretized by positioning the nodes $\x_i \in \Omega$. In the second step, the differential operators $\L$ appearing in the governing equations, are approximated for the given approximation basis. Only then can the system of the governing PDEs be discretized in the third step and finally either assembled into a large sparse system whose solution $\widehat{u}$ is declared as the numerical solution of the governing problem, or solved explicitly by iteration in the fourth step.

To obtain the hybrid scattered-regular domain discretization, we first fill the entire domain with regular nodes. A portion of this regular nodes is then removed in areas where a scattered node placement is desired, i.e., close to the irregular boundaries. Finally, the voids are filled with a dedicated node positioning algorithm~\cite{slak2019generation} that supports variable nodal density and allows us to refine the solution near irregularities in the domain. This approach is rather naive but sufficient for demonstration purposes. A special hybrid fill algorithm is left for future work.

\begin{figure}
	\centering
	\includegraphics[width=\textwidth]{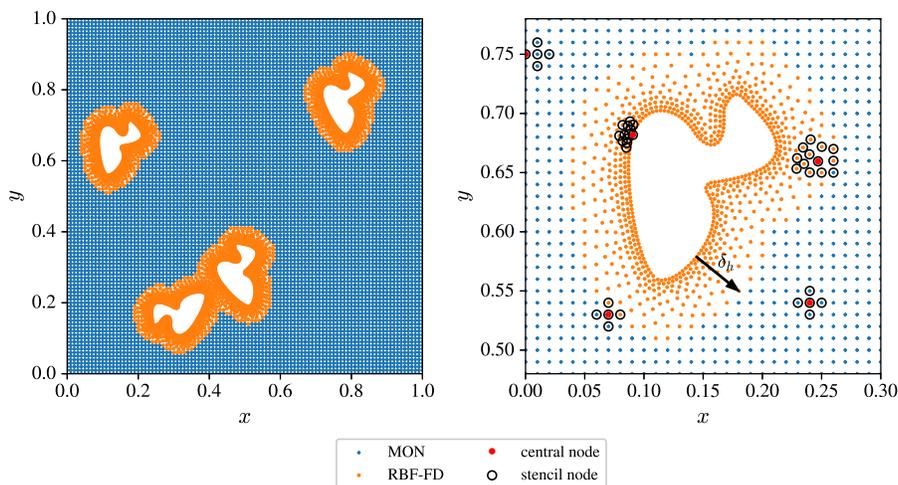}
	\caption{Irregular domain discretization example {\it(left)} and spatial distribution of approximation methods along with corresponding example stencils {\it(right)}.}
	\label{fig:discretization_sample}
\end{figure}

An example of an $h$-refined domain discretization is shown in Fig.~\ref{fig:discretization_sample}. It is worth noting that the width of the scattered node layer $\delta _h$ is non-trivial and affects both the stability of the solution procedure and the accuracy of the numerical solution. Although we provide a superficial analysis of $\delta _h$ variation in Sections~\ref{sec:dvd}~and~\ref{sec:irregular}, further work is warranted.

%To discretize the domain with scattered nodes, we use a dedicated node positioning algorithm~\cite{slak2019generation} also available as a stand-alone implementation in Medusa library~\cite{slak2021medusa}. The chosen node positioning algorithm also conveniently supports a variable nodal density distributions allowing us to employ $h$-refined solution procedures in Sect.~\ref{sec:irregular}. An example of $h$-refined domain discretization using both regular and scattered nodes is shown in Fig.~\ref{fig:discretization_sample}.

After the computational nodes $\x_i\in \Omega$ are obtained, the differential operators $\L$ can be locally approximated in point $\x_c$ over a set of $n$ neighbouring nodes (stencil) $\left \{ \x_i \right \}_{i=1}^n = \mathcal{N}$,  using the following expression
\begin{equation}
	\label{eq:ansatz}
	(\L u)(\x _c) \approx \sum_{i=1}^n w_iu(\x _i).
\end{equation}
The approximation~\eqref{eq:ansatz} holds for an arbitrary function $u$ and yet to be determined weights $\w$. To determine the weights, the equality of approximation~\eqref{eq:ansatz} is enforced for a chosen set of basis functions. Here we will use two variants
\begin{enumerate}[label=(\roman*)]
	\item a set of Polyharmonic splines (PHS) augmented with monomials to ensure convergent behaviour~\cite{flyer2016role,bayona2017role}, effectively resulting in a popular \emph{radial basis function-generated finite differences} (RBF-FD) approximation method~\cite{tolstykh2003using}.
	\item a set of monomials centred at the stencil nodes that we will refer to as (MON)~\cite{kosec2018local}.
\end{enumerate}

We use the least expensive MON with $2d + 1$ monomial basis functions\footnote{In 2D, the 5 basis functions are $\{1, x, y, x^2, y^2\}$. The $xy$ term is not required for regularly placed nodes and its omission allows us to use the smaller and completely symmetric 5-node stencil.} and the same number of support nodes in each approximation stencil. This setup is fast, but only stable on regular nodes~\cite{kosec2018local,slak2019refinedCauchy}. For the RBF-FD part, we also resort to the minimal configuration required for 2nd-order operators, i.e., 3rd-order PHS augmented with all monomials up to the 2nd-order ($m=2$). According to the standard recommendations~\cite{bayona2017role}, this requires a stencil size of $n = 2\binom{m+d}{m}$.

Note the significant difference between stencil sizes --- 5 vs.~12 nodes in 2D --- that only increases in higher dimensions (7 vs.~30 in 3D). This results both in faster computation of the weights $\w$ --- an $\mathcal{O}(N^3)$\footnote{$N_{\mathrm{RBF-FD}} \sim 3 N_{\mathrm{MON}}$ due to the larger stencil size and the extra PHS in the approximation basis.} operation performed only once for each stencil --- and in faster evaluation for the $\mathcal{O}(n)$ explicit operator approximation~\eqref{eq:ansatz} performed many times during the explicit time stepping. Therefore, a spatially varying node regularity can have desirable consequences on the discretization-related errors and computational efficiency of the solution procedure.

\subsection{Computational stability}

By enforcing the equality of approximation~\eqref{eq:ansatz}, we obtain a linear system $\mathbf{M}\b w = \b \ell$. Solving the system provides us with the approximation weights $\b w$, but the stability of such procedure can be uncertain and is usually estimated via the condition number $\kappa(\mathbf{M}) = \left \| \mathbf{M} \right \| \left \| \mathbf{M}^{-1} \right \|$ of matrix $\mathbf{M}$, where $\left \| \cdot \right \|$ denotes the $L^2$ norm.

\begin{figure}
	\centering
	\includegraphics[width=\textwidth]{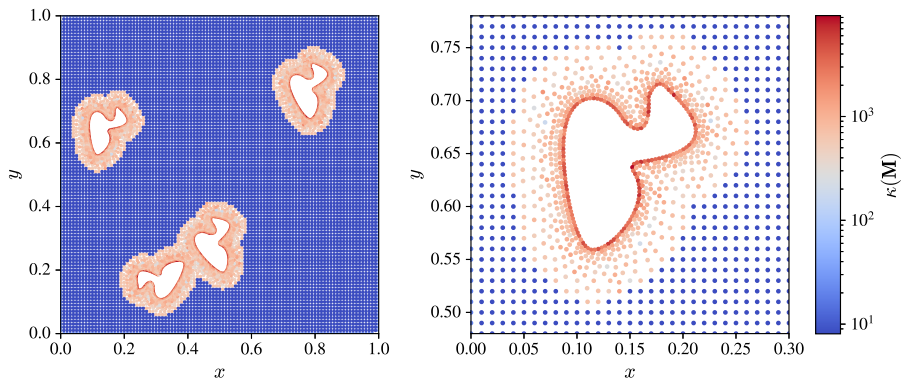}
	\caption{Condition numbers $\kappa(\mathbf{M})$ for the Laplacian operator: entire computational domain {\it(left)} and a zoomed-in section around the irregularly shaped obstacle {\it(right)}.}
	\label{fig:condition_numbers}
\end{figure}

A spatial distribution of condition numbers is shown in Fig.~\ref{fig:condition_numbers}. It can be observed that the RBF-FD approximation method generally results in higher condition numbers than the MON approach. This could be due to the fact that the matrices $\mathbf{M}$ for the RBF-FD part are significantly larger and based on scattered nodes. Nevertheless, it is important to note that the transition from regular to scattered nodes does not appear to affect the conditionality of the matrices.

\subsection{Implementation details}
The entire solution procedure employing the hybrid scattered-regular method is implemented in C++. The projects implementation\footnote{Source code is available at \url{gitlab.com/e62Lab/public/2023_cp_iccs_hybrid_nodes} under tag \emph{v1.1}.} is strongly dependent on our in-house developed meshless C++ framework \emph{Medusa library}~\cite{slak2021medusa} supporting all building blocks of the solution procedure, i.e., differential operator approximations, node positioning algorithms, etc.

We used \texttt{g++ 11.3.0 for Linux} to compile the code with \texttt{-O3 -DNDEBUG} flags on \texttt{Intel(R) Xeon(R) CPU E5520} computer. To improve the timing accuracy we run the otherwise parallel code in a single thread with the CPU frequency fixed at 2.27 GHz, disabled boost functionality and assured CPU affinity using the \texttt{taskset} command.
Post-processing was done using Python 3.10.6 and Jupyter notebooks, also available in the provided git repository.

\section{Governing problem}
\label{sec:examples}
To objectively assess the advantages of the hybrid method, we focus on the natural convection problem that is governed by a system of three PDEs that describe the continuity of mass, the conservation of momentum and the transfer of heat
\begin{align}
	\div \vec{v}                                    & = 0, \label{eq:physics1}                                                                            \\
	\pdv{\vec{v}}{t} + \vec{v} \cdot \grad{\vec{v}} & = -\grad p +\div(\text{Pr} \grad \vec{v}) -\text{Ra}\text{Pr} \vec{g} T_\Delta, \label{eq:physics2} \\
	\pdv{T}{t} + \vec{v} \cdot \grad{T}             & = \div( \grad T), \label{eq:physics3}
\end{align}
where a dimensionless nomenclature using Rayleigh (Ra) and Prandtl (Pr) numbers is used~\cite{de1983natural,kosec2008solution}.

%The dynamics is driven by buoyancy force, which is relatively weak and ensures that the maximum velocity remains well below the speed of sound. As a consequence, the fluid is modelled as incompressible~\cite{anderson2010fundamentals} and allows us to reduce the continuity equation to~\eqref{eq:physics1}. The fluid motion is described by the Navier-Stokes equation~\eqref{eq:physics2}, modified with an additional force term to account for the buoyancy caused by the thermal expansion. Based on the assumption that the acceleration of a fluid driven by natural convection remains insignificant compared to gravity, this force is approximated by the Boussinesq approximation~\cite{Tritton1988}. 
The temporal discretization of the governing equations is solved with the explicit Euler time stepping where we first update the velocity using the previous step temperature field in the Boussinesq term~\cite{Tritton1988}. The pressure-velocity coupling is performed using the Chorin's projection method~\cite{chorin1968numerical} under the premise that the pressure term of the Navier-Stokes equation can be treated separately from other forces and used to impose the incompressibility condition.
%For a detailed description of the temporal discretization and pressure correction consideration the reader is directed to~\cite{rot}. 
The time step is a function of internodal spacing $h$, and is defined as $\mathrm{d}t = 0.1\frac{h^2}{2}$ to assure stability.

% Therefore, we express the temporal derivative with a generalized finite difference approximation to allow for explicit stepping and gather all non-pressure forces in $\mathcal{F}$
% \begin{equation}
% 	\frac{\rho}{\Delta t} \left( \b v(t + \Delta t) - \b v(t) \right) =
% 	-\grad p + \mathcal{F}(\b v(t)).
% \end{equation}
% A single time step is split into two parts, the first resulting in an intermediate velocity that neglects the pressure gradient
% \begin{equation}
% 	\b v^* = \b v(t) + \frac{\Delta t}{\rho} \mathcal{F}(\b v(t))
% \end{equation}
% and the second applying a pressure correction to the intermediate step
% \begin{equation}
% 	\label{eq:chorin2}
% 	\b v(t + \Delta t) = \b v^* -\frac{\Delta t}{\rho} \grad p.
% \end{equation}
% We refer to the pressure term as correction because it does not represent the full hydrodynamic pressure. The correction pressure is calculated by solving the Poisson equation obtained by applying divergence to equation~\ref{eq:chorin2} and imposing equation~\ref{eq:physics1} for the new velocity to ensure incompressibility
% \begin{equation}
% 	\cancelto{0}{\div \b v(t + \Delta t)} - \div \b v^*
% 	= - \frac{\Delta t}{\rho}\laplacian p.
% \end{equation}
% The new velocity is then used in the advection part of the heat transfer equation, which we solve to calculate the new temperature field.

\section{Numerical results}
The governing problem presented in Section~\ref{sec:examples} is solved on different geometries employing (i) MON, (ii) RBF-FD and (iii) their spatially-varying combination. The performance of each approach is evaluated in terms of accuracy of the numerical solution and execution times. Unless otherwise specified, the MON method is employed using the monomial approximation basis omitting the mixed terms, while the RBF-FD approximation basis consists of Polyharmonic splines or order $k=3$ augmented with monomials up to and including order $m=2$.

\subsection{The de Vahl Davis problem}
\label{sec:dvd}
First, we solve the standard  de~Vahl~Davis benchmark problem~\cite{de1983natural}. The main purpose of solving this problem is to establish confidence in the presented solution procedure and to shed some light on the behaviour of considered approximation methods, the stability of the solution procedure and finally on the computational efficiency. Furthermore, the de~Vahl~Davis problem was chosen as the basic test case, because the regularity of the domain shape allows us to efficiently discretize it using exclusively scattered or regular nodes and compare the solutions to that obtained with the hybrid scattered-regular discretization.

\begin{figure}
	\centering
	\includegraphics[width=0.5\textwidth]{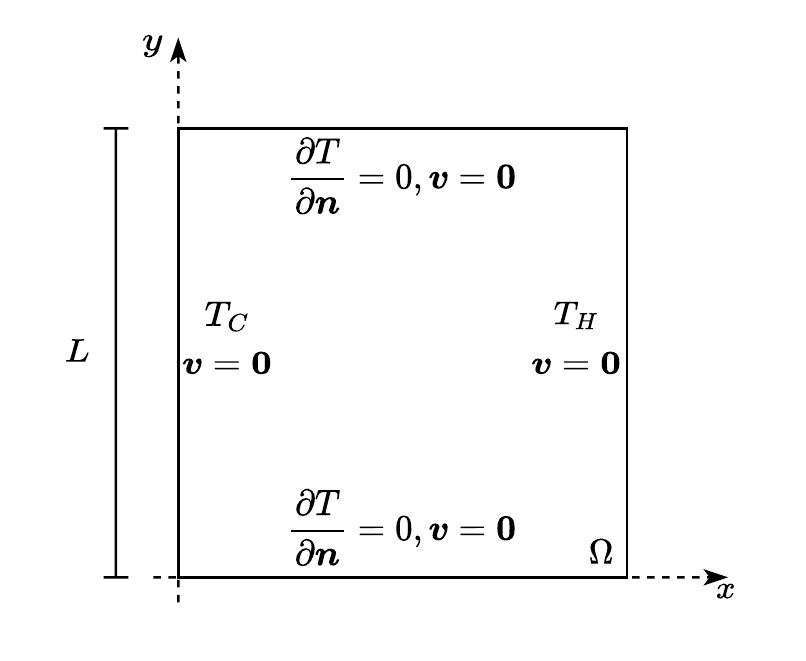}
	\includegraphics[width=0.38\textwidth]{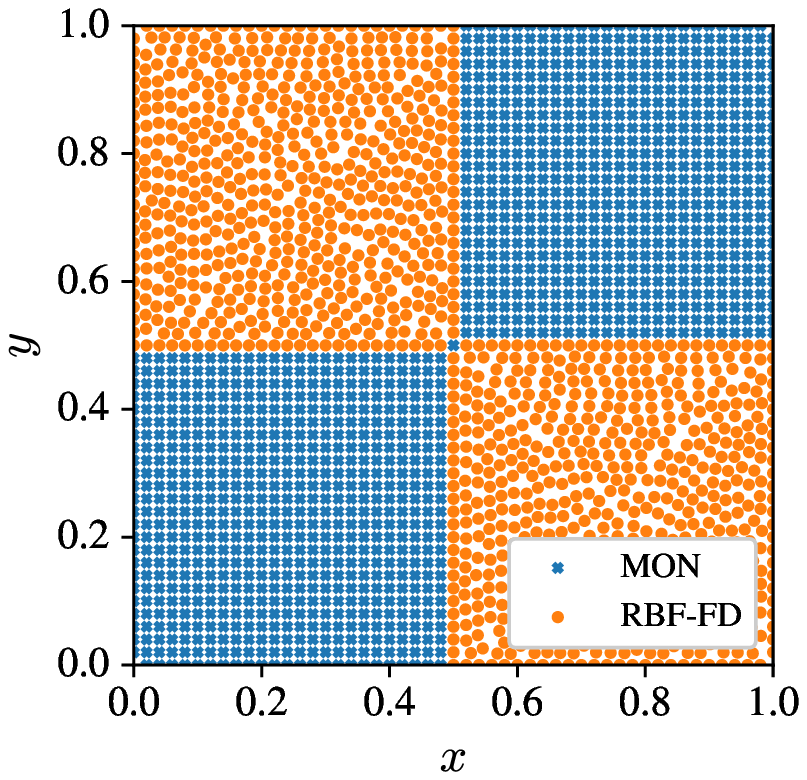}
	\caption{The de Vahl Davis sketch {\it(left)} and example hybrid scattered-regular domain discretization {\it(right)}.}
	\label{fig:dvd_domain}
\end{figure}

For a schematic representation of the problem, see Fig.~\ref{fig:dvd_domain} {\it(left)}. The domain is a unit box $\Omega = \left [ 0, 1 \right ] \times \left [ 0, 1 \right ]$, where the left wall is kept at a constant temperature $T_C = -0.5$, while the right wall is kept at a higher constant temperature $T_H = 0.5$. The upper and lower boundaries are insulated, and no-slip boundary condition for velocity is imposed on all walls. Both the velocity and temperature fields are initially set to zero.

To test the performance of the proposed hybrid scattered-regular approximation method, we divide the domain $\Omega$ into quarters, where each quarter is discretized using either scattered or regular nodes -- see Fig.~\ref{fig:dvd_domain} {\it(right)} for clarity.

\begin{figure}
	\centering
	\includegraphics[width=\textwidth]{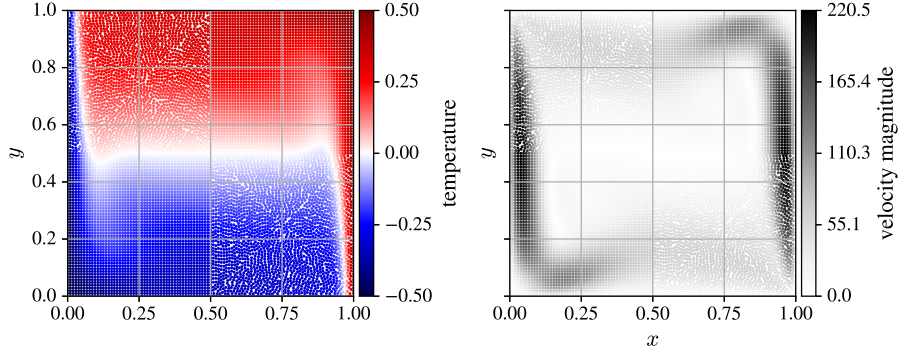}
	\caption{Example solution at the stationary state. Temperature field {\it(left)} and velocity magnitude {\it(right)}.}
	\label{fig:dvd_solution}
\end{figure}

An example solution for $\text{Ra} = 10^6$ and $\text{Pr}=0.71$ at a dimensionless time $t=0.15$ with approximately $N=15\,800$ discretization nodes is shown in Fig.~\ref{fig:dvd_solution}.

\begin{figure}
	\centering
	\includegraphics[width=0.7\textwidth]{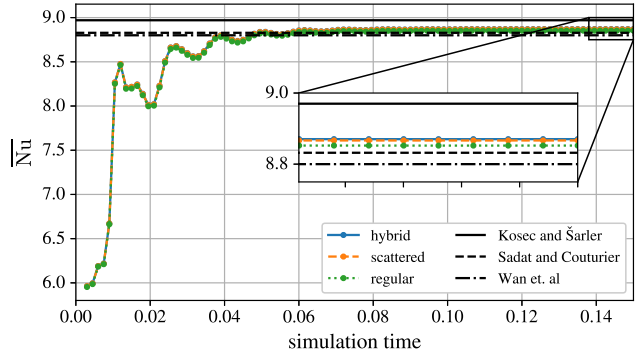}
	\caption{Time evolution of the average Nusselt number along the cold edge calculated with the densest considered discretization. Three reference results Kosec and Šarler~\cite{kosec2008solution}, Sadat and Couturier~\cite{sadat} and Wan et.~al.~\cite{wan} are also added.}
	\label{fig:dvd_nusselt}
\end{figure}

%Due to the circular flow caused by the natural convection, the heat transfer between the heated and cooled walls significantly increases compared to conduction alone.
We use the Nusselt number --- the ratio between convective and conductive heat transfer --- to determine when a steady state has been reached and as a convenient scalar value for comparison with reference solutions. In the following analyses, the average Nusselt number ($\overline{\mathrm{Nu}}$) is calculated as the average of the Nusselt values at the cold wall nodes
\begin{equation}
	\mathrm{Nu} = \frac{L}{T_H-T_C}\abs{\pdv{T}{\boldsymbol n}}_{x=0}.
\end{equation}
Its evolution over time is shown in Fig.~\ref{fig:dvd_nusselt}. In addition, three reference results are also added to the figure. We are pleased to see that our results are in good agreement with the reference solutions from the literature.

Moreover, Fig.~\ref{fig:dvd_nusselt} also shows the time evolution of the average Nusselt number value for cases where the entire domain is discretized using either scattered or regular nodes. We find that all --- hybrid, purely scattered and purely regular domain discretizations --- yield results in good agreement with the references. More importantly, the hybrid method shows significantly shorter computational time (about 50\%) than that required by the scattered discretization employing RBF-FD, as can be seen in Tab.~\ref{tab:nusselt} for the densest considered discretization with $h = 0.00398$.

\begin{table}
	\centering
	\renewcommand{\arraystretch}{1.1}
	\begin{tabular}{cccc}
		Approximation                                    & $\overline{\text{Nu}}$ & execution time [h] & N       \\ \hline \hline
		scattered                                        & 8.867                  & 6.23               & 55\,477 \\
		regular                                          & 8.852                  & 2.42               & 64\,005 \\
		hybrid                                           & 8.870                  & 3.11               & 59\,694 \\ \hline
		Kosec and Šarler (2007)~\cite{kosec2008solution} & 8.97                   & /                  & 10201   \\
		Sadat and Couturier (2000)~\cite{sadat}          & 8.828                  & /                  & 22801   \\
		Wan et.~al.~(2001)~\cite{wan}                    & 8.8                    & /                  & 10201   \\
	\end{tabular}
	\caption{Average Nusselt along the cold edge along with execution times and number of discretization nodes.}
	\label{tab:nusselt}
\end{table}

\begin{figure}
	\centering
	\includegraphics[width=\textwidth]{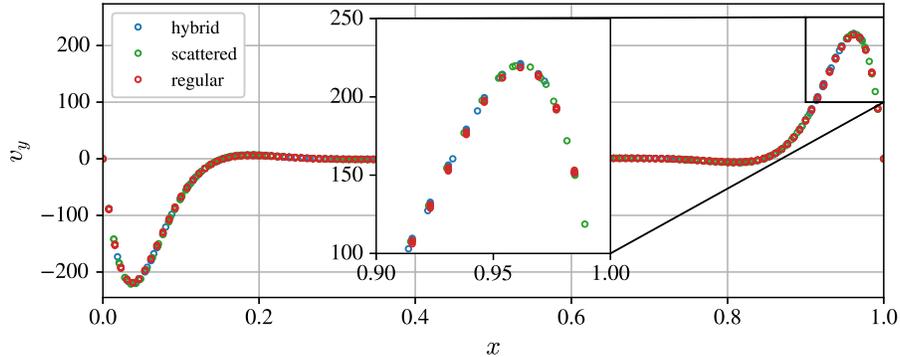}
	\caption{Vertical velocity component values at nodes close to the vertical midpoint of the domain, i.e., $|y - 0.5| \le h$ for purely scattered, purely regular and hybrid discretizations.}
	\label{fig:dvd_sym}
\end{figure}

To further validate the hybrid method, we show in Fig.~\ref{fig:dvd_sym} the vertical component of the velocity field across the section $y=0.5$. It is important to observe that the results for the hybrid, scattered and regular approaches overlap, which means that the resulting velocity fields for the three approaches are indeed comparable.

\begin{figure}
	\centering
	\includegraphics[width=\textwidth]{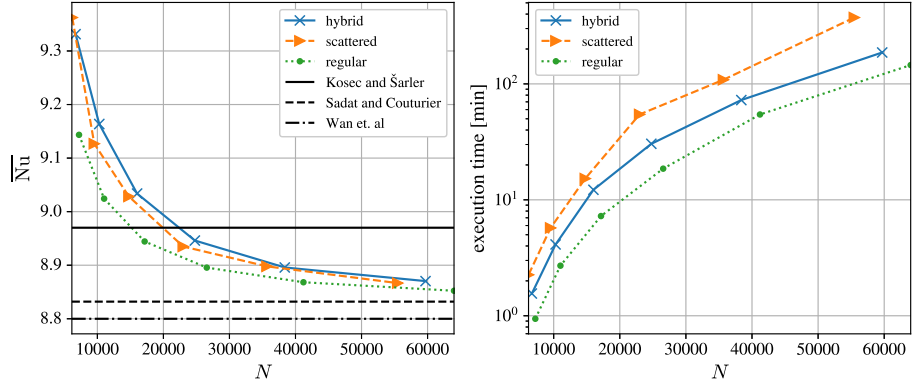}
	\caption{Convergence of average Nusselt number with respect to discretization quality {\it(left)} and corresponding execution times {\it(right)}.}
	\label{fig:dvd_conv}
\end{figure}

As a final remark, we also study the convergence of the average Nusselt number with respect to the number of discretization nodes in Fig.~\ref{fig:dvd_conv} {\it(left)}, where we confirm that all our discretization strategies converge to a similar value that is consistent with the reference values. Moreover, to evaluate the computational efficiency of the hybrid approach, the execution times are shown on the right. Note that the same values for $h$ were used for all discretization strategies and the difference in the total number of nodes is caused by the lower density of scattered nodes at the same internodal distance.

\subsubsection{The effect of the scattered nodes layer width $\delta_h$}

To study the effect of the width of the scattered node layer $\delta_h$, we consider two cases. In both cases, the domain from Fig.~\ref{fig:dvd_domain} is split into two parts at a distance $h\delta_h$ from the origin in the lower left corner. In the first scenario, the split is horizontal, resulting in scattered nodes below the imaginary split and regular nodes above it. In the second scenario, the split is vertical, resulting in scattered nodes to the left of it and regular nodes to the right of it. In both cases, the domain is discretized with purely regular nodes when $h\delta_h = 0$ and with purely scattered nodes when $h\delta_h = L$.

\begin{figure}
	\centering
	\includegraphics[width=\textwidth]{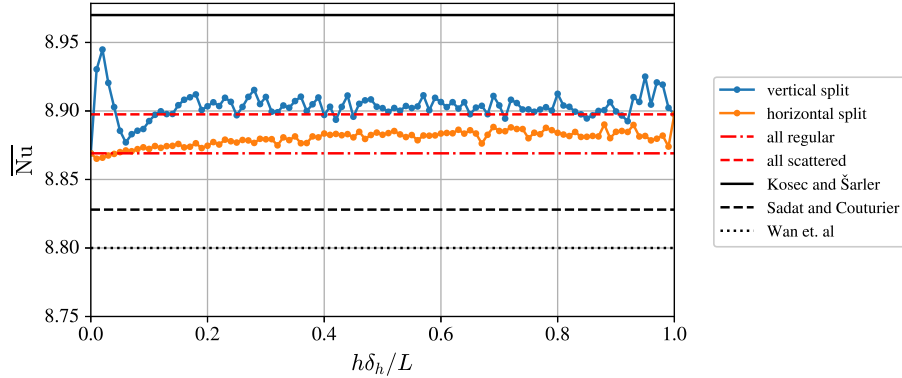}
	\caption{Demonstration of the scattered node layer width ($\delta_h$) effect on the accuracy of the numerical solution.}
	\label{fig:dvd_delta}
\end{figure}

In Fig.~\ref{fig:dvd_delta}, we show how the width of the scattered node layer affects the average Nusselt number in stationary state for approximately $40\,000$ discretization nodes. It is clear that even the smallest values of $\delta_h$ yield satisfying results. However, it is interesting to observe that the accuracy is significantly affected when the boundary between regular and scattered nodes runs across the region with the largest velocity magnitudes, i.e., the first and last couple of vertical split data points in Fig.~\ref{fig:dvd_delta}.

\subsection{Natural convection on irregularly shaped domains}
\label{sec:irregular}
In the previous section we demonstrated that the hybrid scattered-regular approximation method is computationally more efficient than the pure RBF-FD approximation with only minor differences in the resulting fields. However, to truly exploit the advantages of the hybrid method, irregular domains must be studied. Therefore, in this section, the hybrid scattered-regular approach is employed on an irregularly shaped domain. Let the computational domain $\Omega$ be a difference between the two-dimensional unit box $\Omega = \left [ 0, 1 \right ] \times \left [ 0, 1 \right ]$ and 4 randomly positioned duck-shaped obstacles introducing the domain irregularity.

The dynamics of the problem are governed by the same set of equations~(\ref{eq:physics1}-\ref{eq:physics3}) as in the previous section. This time, however, all the boundaries of the box are insulated. The obstacles, on the other hand, are subject to Dirichlet boundary conditions, with half of them at $T_C=0$ and the other half at $T_H=1$. The initial temperature is set to $T_{\text{init}} = 0$.

We have chosen such a problem because it allows us to further explore the advantages of the proposed hybrid scattered-regular discretization. Generally speaking, the duck-shaped obstacles within the computational domain represent an arbitrarily complex shape that requires scattered nodes for accurate description, i.e., reduced discretization-related error. Moreover, by using scattered nodes near the irregularly shaped domain boundaries, we can further improve the local field description in their vicinity by employing a $h$-refined discretization. Specifically, we employ $h$-refinement towards the obstacles with linearly decreasing internodal distance from $h_r=0.01$ (regular nodes) towards $h_s=h_r/3$ (irregular boundary) over a distance of $h_r \delta_h$. The refinement distance and the width of the scattered node layer are the same, except in the case of fully scattered discretization. Such setup effectively resulted in approximately $N=11\,600$ computational nodes ($N_s=3149$ scattered nodes and $N_r = 8507$ regular nodes), as shown in Fig.~\ref{fig:discretization_sample} for a scattered node layer width $\delta_h = 4$. Note that the time step is based on the smallest $h$, i.e., $\mathrm{d}t = 0.1\frac{h_s^2}{2}$.

\begin{figure}
	\centering
	\includegraphics[width=\textwidth]{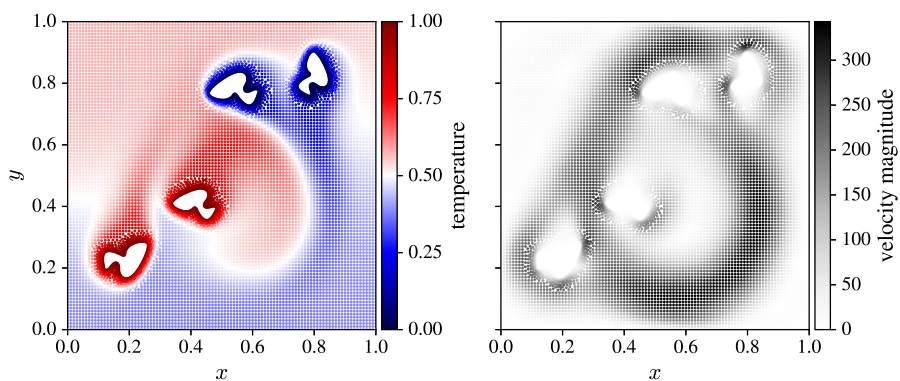}
	\caption{Example solution on irregular domain. Temperature field {\it(left)} and velocity magnitude {\it(right)}.}
	\label{fig:natural_solution}
\end{figure}

Fig.~\ref{fig:natural_solution} shows an example solution for an irregularly shaped domain. The hybrid scattered-regular solution procedure was again able to obtain a reasonable numerical solution.

\begin{figure}
	\centering
	\includegraphics[width=\textwidth]{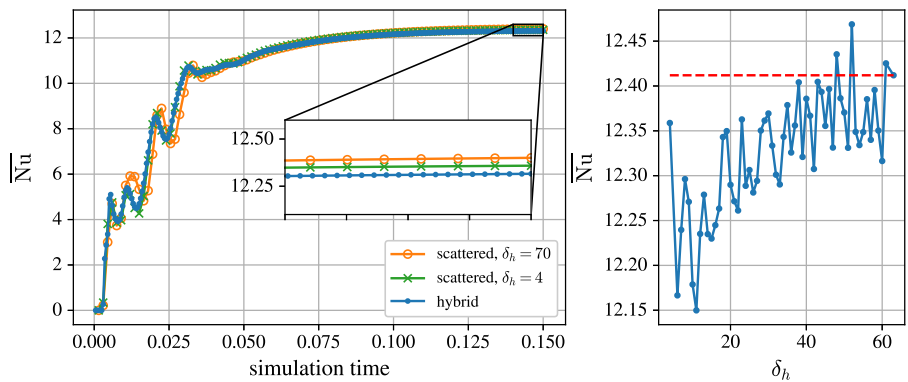}
	\caption{Time evolution of the average Nusselt number calculated on the cold duck-shaped obstacles of an irregularly shaped domain. {\it(left)} Changes in stationary state average Nusselt number as the scattered node layer width $\delta_h$ increases. {\it(right)}}
	\label{fig:irregular_nusselt}
\end{figure}

Furthermore, Fig.~\ref{fig:irregular_nusselt} {\it(left)} shows the average Nusselt number along the cold duck edges where we can observe that a stationary state has been reached. The steady state values for all considered discretizations match closely but it is interesting to note that in the early stage of flow formation, the fully scattered solutions with different refinement distance $\delta_h$ differ significantly more than the hybrid and the fully scattered solutions with the same refinement strategy.

It is perhaps more important to note that the execution times gathered in Tab.~\ref{tab:nusselt_ireg} show that the hybrid method effectively reduces the execution time for approximately 35~\%. The pure regular discretization with MON approximation is omitted from the table because a stable numerical solution could not be obtained.

\begin{table}[h]
	\centering
	\renewcommand{\arraystretch}{1.1}
	\begin{tabular}{cccc}
		Approximation & $\overline{\text{Nu}}$ & execution time [min] & N       \\ \hline \hline
		scattered     & 12.32                  & 46.31                & 10\,534 \\
		hybrid        & 12.36                  & 29.11                & 11\,535 \\
	\end{tabular}
	\caption{Average Nusselt along the cold duck edges along with execution times. Note that all values in the table were obtained for $\delta_h = 4$.}
	\label{tab:nusselt_ireg}
\end{table}

\subsubsection{The effect of the scattered nodes layer width $\delta_h$}

To justify the use of $\delta_h = 4$, we show in Fig.~\ref{fig:irregular_nusselt} {\it(right)} that the average value of the Nusselt number at steady state for different values of $\delta_h$. In the worst case, the difference is $<2 \%$, justifying the use of the computationally cheaper smaller $\delta_h$. Note that in this particular domain setup, $\delta_h > 64$ already yields a purely scattered domain discretization, while the minimum stable value is $\delta_h = 4$. Note also that the general increase of the Nusselt number with respect to the width of the scattered node layers $\delta_h$ may also exhibit other confounding factors. An increase in $\delta_h$ leads to a finer domain discretization due to a more gradual refinement, i.e., a fully scattered discretization using $\delta_h = 70$ results in about 35\,000 discretization points compared to 11\,600 at $\delta_h=4$, while a decrease in $\delta_h$ leads to a more aggressive refinement that could also have a negative effect. This can be supported by observing the difference between the results for the two fully scattered discretizations in Fig.~\ref{fig:irregular_nusselt}.

% \begin{figure}
% 	\centering
% 	\includegraphics[width=0.8\textwidth]{images/duck_nusselt_delta_h.png}
% 	\caption{Effect of the scattered node layer width ($\delta_h$) around the irregular boundaries on the average Nusselt number computed on the cold boundaries.\mr{A dam ta graf side-by-side s prejšnim? Mogoče 2:1 širina? Lahko bi kazala še različne $\delta_h$ pri fully scattered, da bi bilo malo manj filozofiranja pri vpeljavi slike.}}
% 	\label{fig:duck_delta_h}
% \end{figure}

% In Fig.~\ref{fig:dvd_delta} we show how the width of scattered node layer affects the average Nusselt number in stationary state for approximately $40\,000$ discretization nodes. Clearly, even the smallest values of $\delta_h$ yield satisfying results. However, it is interesting to observe that the accuracy is notably affected if the boundary between regular and scattered nodes runs across the area with the largest fluid velocity magnitudes, i.e, the first and last couple of vertical split data points in Fig.~\ref{fig:dvd_delta}.

% \begin{figure}
% 	\centering
% 	\includegraphics[width=\textwidth]{images/dvd_nusselt_delta_h.png}
% 	\caption{Demonstration of the scattered node layer width ($\delta_h$) effect on the accuracy of the numerical solution.}
% 	\label{fig:dvd_delta}
% \end{figure}

%
%
\subsection{Application to three-dimensional irregular domains}
\label{sec:irregular3d}
As a final demonstrative example, we employ the proposed hybrid scattered-regular approximation method on a three-dimensional irregular domain. The computational domain $\Omega$ is a difference between the three-dimensional unit box $\Omega = \left [ 0, 1 \right ] \times \left [ 0, 1 \right ]\times \left [ 0, 1 \right ]$ and 4 randomly positioned and sized spheres introducing the domain irregularity.

The dynamics are governed by the same set of equations (\ref{eq:physics1}-\ref{eq:physics3}) as in the two-dimensional case from Sect.~\ref{sec:irregular}. To improve the quality of the local field description near the irregularly shaped domain boundaries, $h$-refinement is employed with a linearly decreasing internodal distance from $h_r=0.025$ (regular nodes) towards $h_s=h_r/2$ (spherical boundaries). Two spheres were set to a constant temperature $T_C = 0$ and the remaining two to $T_H=1$. The Rayleigh number was set to $10^4$.

\begin{figure}
	\centering
	\includegraphics[width=\textwidth]{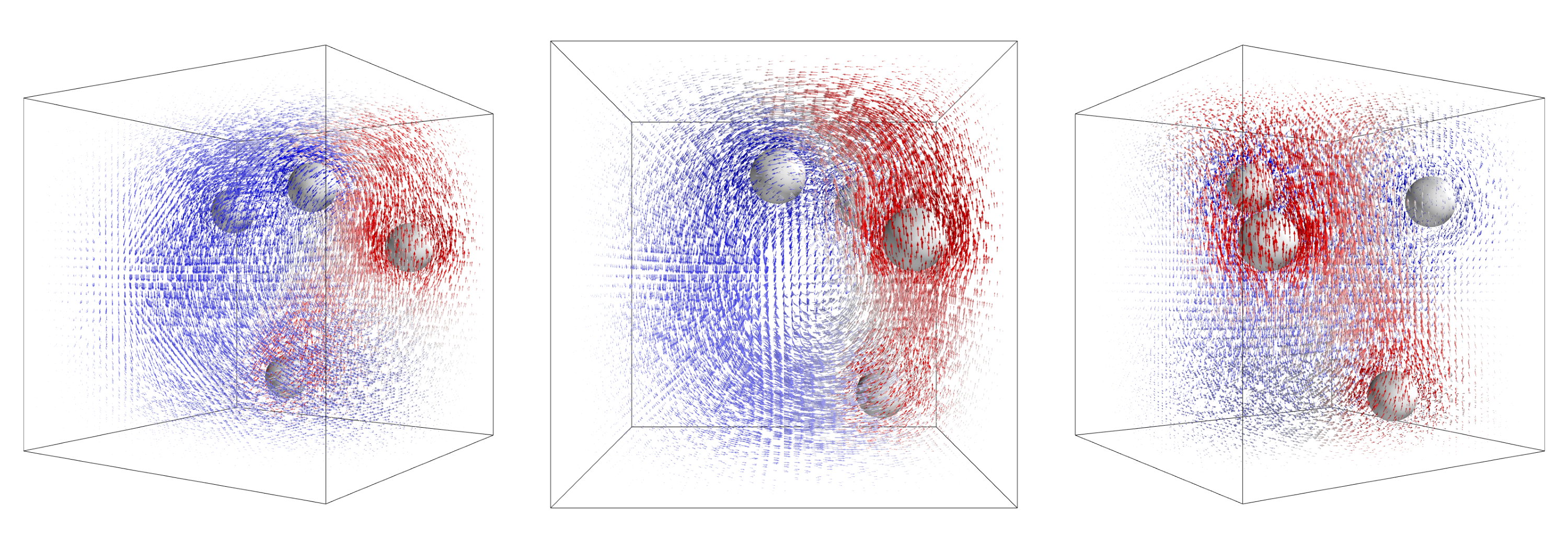}
	\caption{Example solution viewed from three different angles. The arrows show the velocity in computational nodes and are coloured according to the temperature in that node. The values range from dark blue for $T_C$ to dark red for $T_H$. For clarity, only a third of the nodes is visualized.}
	\label{fig:3d_solution}
\end{figure}

Although difficult to visualize, an example solution is shown in Fig.~\ref{fig:3d_solution}. Using the hybrid scattered-regular domain discretization, the solution procedure was again able to obtain a reasonable numerical solution.

%Moreover, the time evolution of the average Nusselt number along the hot spherical boundaries, shown in Fig.~\ref{fig:3d_nusselt}\mr{Glede tega grafa nisva sigurna, ce bi ga sploh kazala.}, suggests that both the hybrid and the purely scattered approaches were able to obtain a steady-state solution, while using a purely regular discretization again failed to yield a numerical solution. 
Note that the scattered method took about 48 hours and the hybrid scattered-regular approximation method took 20 hours to simulate 1 dimensionless time unit with the dimensionless time step $\mathrm{d}t = 7.8125\cdot10^{-6}$ and about 75\,000 computational nodes with $\delta_h = 4$. For clarity, the data is also gathered in Tab.~\ref{tab:nusselt_ireg_3d}. Note that the pure regular discretization with MON approximation is again omitted from the table because a stable numerical solution could not be obtained.

\begin{table}[h]
	\centering
	\renewcommand{\arraystretch}{1.1}
	\begin{tabular}{cccc}
		Approximation & $\overline{\text{Nu}}$ & execution time [h] & N       \\ \hline \hline
		scattered     & 7.36                   & 48.12              & 65\,526 \\
		hybrid        & 6.91                   & 20.54              & 74\,137 \\
	\end{tabular}
	\caption{Average Nusselt along the cold spheres, execution time, and number of computational nodes.}
	\label{tab:nusselt_ireg_3d}
\end{table}

\section{Conclusions}
We proposed a computationally efficient approach to the numerical treatment of problems in which most of the domain can be efficiently discretized with regularly positioned nodes, while scattered nodes are used near irregularly shaped domain boundaries to reduce the discretization-related errors. The computational effectiveness of the spatially-varying approximation method, employing FD-like approximation on regular nodes and RBF-FD on scattered nodes, is demonstrated on a solution to a two-dimensional de Vahl Davis natural convection problem.

We show that the proposed hybrid method, can significantly improve the computational efficiency compared to the pure RBF-FD approximation, while introducing minimal cost to the accuracy of the numerical solution. A convergence analysis from Fig.~\ref{fig:dvd_conv} shows good agreement with the reference de~Vahl~Davis solutions.

In the continuation, the hybrid method is applied to a more general natural convection problem in two- and three-dimensional irregular domains, where the elegant mathematical formulation of the meshless methods is further exposed by introducing $h$-refinement towards the irregularly shaped obstacles. In both cases, the hybrid method successfully obtained the numerical solution and proved to be computationally efficient, with execution time gains nearing 50~\%.

Nevertheless, the scattered node layer width and the aggressiveness of $h$-refinement near the irregularly shaped domain boundaries should be further investigated, as both affect the computational efficiency and stability of the solution procedure. In addition, future work should also include more difficult problems, such as mixed convection problems and a detailed analysis of possible surface effects, e.g.\ scattering, at the transition layer between the scattered and regular domains.
\subsubsection{Acknowledgements}
The authors would like to acknowledge the financial support of Slovenian Research Agency (ARRS) in the framework of the research core funding No.\ P2-0095, the Young Researcher program PR-10468 and research project J2-3048.

\subsubsection{Conflict of interest}
The authors declare that they have no conflict of interest. All the co-authors have confirmed to know the submission of the manuscript by the corresponding author.

\bibliographystyle{splncs04}
\bibliography{references}

\end{document}